\documentclass{amsart}

\usepackage{breqn}
\usepackage{mathtools}
\usepackage{accents}

\begin{document}

\title[curvature]{Curvature on the integers, I}
\bigskip

\def \h{\hat{\ }}
\def \cO{\mathcal O}
\def \ra{\rightarrow}
\def \bZ{{\mathbb Z}}
\def \cP{{\mathcal V}}
\def \cH{{\mathcal H}}
\def \cB{{\mathcal B}}
\def \d{\delta}
\def \cC{{\mathcal C}}

\newtheorem{THM}{{\!}}[section]
\newtheorem{THMX}{{\!}}
\renewcommand{\theTHMX}{}
\newtheorem{theorem}{Theorem}[section]
\newtheorem{corollary}[theorem]{Corollary}
\newtheorem{lemma}[theorem]{Lemma}
\newtheorem{proposition}[theorem]{Proposition}
\theoremstyle{definition}
\newtheorem{definition}[theorem]{Definition}
\theoremstyle{remark}
\newtheorem{remark}[theorem]{Remark}
\newtheorem{example}[theorem]{\bf Example}
\numberwithin{equation}{section}
\author{Malik Barrett and Alexandru Buium}
\address{Department of Mathematics and Statistics\\University of New Mexico \\ Albuquerque, NM 87131, USA}
\email{buium@math.unm.edu, malikb@gmail.com}

\begin{abstract} 
Starting with a symmetric/antisymmetric matrix with integer coefficients 
(which we view as an analogue of a metric/form on a principal bundle over the ``manifold" $Spec\ \bZ$) we introduce arithmetic analogues of Chern connections and their curvature (in which usual partial derivative operators acting on functions are replaced by Fermat quotient operators acting on integer numbers); curvature is introduced via the method of ``analytic continuation between primes" \cite{laplace}. We prove 
various  non-vanishing, respectively vanishing results for curvature; morally, $Spec\ \bZ$ will appear as ``intrinsically curved." Along with \cite{adel1, adel2, adel3}, this theory can be viewed as taking first steps  in developing   a ``differential geometry of $Spec\ \bZ$". 
\end{abstract}

\maketitle


\section{Introduction}

 In the present paper we will view $\bZ$ as an analogue of a ring of ``algebraic/analytic" functions on an infinite dimensional  manifold in which the various ``directions" correspond to the primes; $C^{\infty}$ objects in geometry will then roughly correspond to adelic objects in arithmetic. In the spirit of \cite{char, book, laplace}, the partial derivative operators, acting on functions, will be replaced by  Fermat quotient type  operators (called {\it $p$-derivations}), acting on numbers.  We will then want to develop an arithmetic analogue  of connections and curvature on  $Spec\ \bZ$. 
In our paper the analogue of $C^{\infty}$ connections on $Spec\ \bZ$ will be families $(\d_p)$ indexed by rational primes $p$ (belonging to  a  fixed, possibly infinite, set ${\mathcal V}$) where for each $p$, $\d_p$ is a $p$-derivation on the $p$-adic completion  of the ring $\cO(GL_n)$ (where $GL_n$ is the general linear group scheme over some ring of cyclotomic integers).  Families $(\d_p)$ as above will be referred to as {\it adelic connections}.   
  One immediately encounters, however, the following difficulty: the
operators $\d_p$, for various $p$'s,  do not act a priori on the same ring and, hence,  one cannot directly consider their commutator; consequently the notion of curvature, which should correspond to such a family of commutators,  seems problematic. We will overcome this problem in two ways: the first, to be explained in the present paper,  is by implementing the technique of {\it analytic continuation between primes} introduced in \cite{laplace}; the second way, to be explained in a sequel to this paper \cite{curvature2},  is by 
 {\it algebraization of Frobenius lifts by correspondences}, a method that  has a 
  ``birational/motivic" flavor. In both paradigms a notion of curvature  can then be developed. Our main results, in this paper and its sequel \cite{curvature2},  are theorems about the vanishing/non-vanishing of the curvatures of some remarkable  adelic connections, referred to as {\it Chern connections}, arising from \cite{adel2}; these Chern connections  are naturally attached to symmetric/antisymmetric matrices (equivalently to outer involutions of $GL_n$ defining various forms of the classical groups $SO_n$ and $Sp_n$), and are an arithmetic analogue of  Chern connections \cite{GH, kobayashi} attached to metrics/forms in differential geometry. 
    
   It is worth pointing out here that what we call {\it curvature} in our  context is entirely  different from what is called {\it $p$-curvature} in the arithmetic theory of differential equations that has been developed around the Grothendieck conjecture (cf., e.g., \cite{katz}); indeed our curvature here is about the ``$p$-differentiation" of numbers with respect to primes $p$ (in other words it is about $d/dp$) whereas the theory in \cite{katz} (and related papers) is about usual  differentiation $d/dt$ with respect to a variable $t$ of power series in $t$ with arithmetically interesting coefficients. In spite of these differences the two types of curvatures could interact; a model for such an interaction between $d/dp$ and $d/dt$ is in the papers \cite{pde,pdemod}. 
   
   Along with \cite{adel1, adel2, adel3}, the present paper and its sequel \cite{curvature2}   can be viewed as taking first steps  in developing   a ``differential geometry of $Spec\ \bZ$."
   The theory can be further developed by considering various (analogues of) contractions of curvature such as Ricci, mean, and scalar curvature; this will be done elsewhere. 
Another direction is the arithmetic analogue of the  symplectic/Hamiltonian story introduced in \cite{BYM}.
   The whole project is closely related to  Borger's viewpoint  on the field with one element \cite{borgerf1}; cf. also the Introduction to \cite{book}.

Our paper is organized as follows. In section 2 we introduce our main arithmetic concepts preceded by  a motivational discussion of the corresponding classical differential geometric concepts. In section 3 we state our main results. In section 4 we prove our results. In section 5 we add some remarks on the case when ${\mathcal V}$ consists of  {\it one prime $p$ only}; this case becomes, through the consideration of the $p$-derivations $\d_p$ and  ``$\d_{\overline{p}}$," analogous to the consideration of  connections  of vector (or principal) bundles on complex curves.

\medskip

{\bf Acknowledgment} 
The second author would like to acknowledge inspiring discussions with J. Borger, C. Boyer, Yu. I. Manin, and D. Vassilev. Also, during the preparation of this paper,  the second author was 
 partially supported by the Simons Foundation
(award 311773),   the Institut des Hautes Etudes Scientifiques in Bures sur Yvette,    the Romanian National Authority
for Scientific Research (CNCS - UEFISCDI, 
PN-II-ID-PCE-2012-4-0201), and  the Max-Planck-Institut f\"{u}r Mathematik in Bonn.

\section{Main concepts}

We start with a discussion of some classical concepts in differential geometry. We then proceed to introduce their arithmetic analogues. 

\subsection{Classical differential geometry}
As a motivation for our arithmetic theory we ``recall" some basic concepts from classical differential geometry. From a logical viewpoint this subsection will not play any role later. However our presentation of this classical material is not completely standard and the twist on the standard picture is meant to show the way to the arithmetic case.

\subsubsection{Connections in principal bundles}
Let  $A$ be a ring (commutative, with identity) and let $\d_1,...,\d_m$ be commuting derivations 
 \begin{equation}
 \label{deltai}
 \d_i:A\ra A.
 \end{equation}
 Let $B$ be an $A$-algebra. By a {\it connection} on $B$ (equivalently,  on the scheme $Spec \ B$) we will understand an $m$-tuple  $\d=(\d_i)$ where for each $i$,  $\d_i:B\ra B$ is a derivation  lifting the corresponding derivation  on $A$. We are especially interested in the following special case.
 Denote by $G=GL_n=Spec\ B$, with $B:=A[x,\det(x)^{-1}]$, the general linear group scheme over $A$, where $x=(x_{kl})$ is an $n\times n$ matrix of indeterminates.
  So $G(A)=GL_n(A)$ is the group of all invertible $n\times n$ matrices with entries in $A$ and 
 $\cO(G)=A[x,\det(x)^{-1}]$. Denote   by 
 ${\mathfrak g}=Spec\ A[x]$ so
 ${\mathfrak g}(A)$ the additive group  of all $n\times n$ matrices with entries in $A$; we view it both as associative algebra and as Lie algebra over $A$. 
 By the above a {\it connection}  on $G$ is an 
  $m$-tuple $\d=(\d_i)$, where 
 \begin{equation}
 \label{godot}
 \d_i:\cO(G)\ra \cO(G)\end{equation}
  are derivations, extending the corresponding derivations on $A$.
  A connection $\d=(\d_i)$ on $G$ will be called {\it right invariant} (or simply {\it invariant}) if
  for all $i$,  $\d_i x=A_i x$, where $A_i\in {\mathfrak g}(A)$.  Here 
  $\d_i x$ is the matrix with entries $\d_i x_{kl}$. More generally, in this paper, 
   if $u_{kl}$ are the entries of a matrix $u$
  with coefficients in a ring and $f$ is a map of that ring into itself then we write $f(u)$ for the matrix $(f(u_{kl}))$ with entries $f(u)_{kl}=f(u_{kl})$.

\subsubsection{Curvature, $3$-curvature, $(1,1)$-curvature} For any   connection $\d=(\d_i)$ on $G$ one
defines the {\it curvature} of $\d$ as the family $(\varphi_{ij})$ where $\varphi_{ij}$ are  
   the $A$-derivations
 \begin{equation}
\label{theta}
\varphi_{ij}:=[\d_i,\d_j]:=\d_i\d_j-\d_j\d_i:\cO(G)\ra \cO(G).\end{equation}
If $\d$ is invariant, with $\d_i x=A_ix$, $A_i\in {\mathfrak g}(A)$, then
\begin{equation}
\label{Fij}
\varphi_{ij}(x)=F_{ij}x,\ \ \ F_{ij}:=\d_iA_j-\d_jA_i-[A_i,A_j]\in {\mathfrak g}(A),\end{equation}
where $[A_i,A_j]=A_iA_j-A_jA_i$ is the usual commutator of matrices.
The matrix $F=(F_{ij})$ is also referred to, in this case, as  the {\it curvature} of the invariant connection $\d=(\d_i)$.  

For a connection $\d=(\d_i)$ one can define the {\it $3$-curvature} of $\d$ as the family $(\varphi_{ijk})$ of $A$-derivations
\begin{equation}
\label{yangmillssymbol}
\varphi_{ijk}:=[\d_i,[\d_j,\d_k]]=[\d_i,\varphi_{jk}]:\cO(G)\ra \cO(G).
\end{equation}
For $\d$ invariant as above, we have
$\varphi_{ijk}(x)=(\d_i F_{jk}-[A_i,F_{jk}])x$. In general, $\varphi_{ijk}$ is antisymmetric in $jk$
 and satisfies the {\it Bianchi identity}:
$\sum_{(ijk)}\varphi_{ijk}=0$,
where $(ijk)$ indicates summation over the cyclic permutations of $i,j,k$. The $3$-curvature plays a role in the Yang-Mills formalism \cite{connesdubois}.

The above discussion  has a ``$(1,1)$ analogue" as follows.
Assume 
  that we have two invariant connections on $G$,
 $$\d=(\d_i)_{1\leq i \leq m},\ \ \ 
\bar{\d}=(\bar{\d}_i)_{1\leq i\leq m}$$
for two $m$-tuples of derivations on $A$ 
such that 
$\{\d_1,...,\d_m,\bar{\d}_1,...,\bar{\d}_m\}$ are pairwise commuting on $A$; we then say that $\overline{\d}$  {\it commutes} with $\d$ on $A$. We
  write $\d_{\bar{i}}:=\bar{\d}_i$. 
 We have at our disposal the curvature $(F_{ij})$  of $\d$
and also the  curvature $(F_{\overline{i}\overline{j}})$ of $\overline{\d}$.
 On the other hand we can consider the commutators
$$\varphi_{i\bar{j}}:=[\d_i,\d_{\bar{j}}]:\cO(G)\ra \cO(G).$$ If $\d_i x=A_ix$, $\d_{\bar{i}}x=A_{\bar{i}}x$, then
$$\varphi_{i\bar{j}}(x)=F_{i\bar{j}}x,\ \ \ F_{i\bar{j}}:=\d_iA_{\bar{j}}-\d_{\bar{j}}A_i-[A_i,A_{\bar{j}}].$$
We call $(\varphi_{i\bar{j}})$, or $(F_{i\bar{j}})$,  the {\it $(1,1)$-curvature} of $\d$ with respect to $\overline{\d}$. 

\subsubsection{Chern connections} 
Let $q\in G(A)=GL_n(A)$ be a symmetric, respectively antisymmetric matrix; hence $q^t=\pm q$, where the $t$ superscript means transpose. Consider the  maps
$\cH_q:G\ra G$, $\cB_q:G\times G\ra G$ defined by
$$\cH_q(x)=x^tqx,\ \ \ \cB_q(x,y)=x^tqy.$$
Consider the  {\it trivial connection} $\d_0=(\d_{0i})$ on $G$ defined by $\d_{0i}x=0$. Then there is a unique  invariant connection $\d=(\d_i)$ on $G$ (which we refer to as the {\it Chern connection} on $G$ attached to $q$) such that the following diagrams are commutative:
 \begin{equation}
 \label{got}
 \begin{array}{ccc}
 \cO(G) & \stackrel{\d_i}{\longleftarrow} & \cO(G)\\
 \cH_q \uparrow &\ &\uparrow \cH_q\\
 \cO(G) & \stackrel{\d_{0i}}{\longleftarrow} & \cO(G)\end{array}\ \ \ \ \ 
 \begin{array}{rcl}
 \cO(G) & \stackrel{\d_i\otimes 1+1\otimes \d_{0i}}{\longleftarrow} & \cO(G)\otimes_A \cO(G)\\
 \d_{0i}\otimes 1+1\otimes \d_i \uparrow & \  & \uparrow \cB_q\\
 \cO(G)\otimes_A \cO(G) & \stackrel{\cB_q}{\longleftarrow} & \cO(G)\end{array}
 \end{equation}
 In the above diagram we still denoted by ${\mathcal H}_q$ and ${\mathcal B}_q$ the ring homomorphisms induced by ${\mathcal H}_q$ and ${\mathcal B}_q$ respectively. 
 To explain this and subsequent remarks it is convenient to switch to some classical notation as follows.
 Let  $\d_ix=A_i x$, set 
 $\Gamma_i :=  -A_i^t$, let $\Gamma_{ij}^k$ be the   $(j,k)$-entry of $\Gamma_i$
(the {\it Cristoffel symbols}), and  set $\Gamma_{ijk}  := \Gamma_{ij}^lq_{lk}$ (Einstein summation). 
 Then the commutativity of the left diagram in \ref{got} is equivalent to
 \begin{equation}
\label{horizontal}
\d_i q_{jk}=\Gamma_{ijk}\pm \Gamma_{ikj};
\end{equation}
classically one refers to this condition as the {\it compatibility of the connection with the metric}.
On the other hand the commutativity of the right diagram in \ref{got} is equivalent to
\begin{equation}
\label{symmetric}
\Gamma_{ijk}=\pm\Gamma_{ikj}.
\end{equation}
So, given $q$,  equations  \ref{horizontal} and \ref{symmetric} have a unique solution $\Gamma_{ijk}$, given by
\begin{equation}
\label{solution}
\Gamma_{ijk}=\frac{1}{2}\d_i q_{jk},
\end{equation}
equivalently by
\begin{equation}
\label{sol}
A_i=-\frac{1}{2} q^{-1}\d_i q.
\end{equation}
Our Chern connection $\d$ attached to $q$ is an analogue of hermitian Chern connections on hermitian vector bundles, cf. \cite{GH}, p. 73., and also of the Duistermaat connections  \cite{dui} in the real case. These analogies  are explained in detail in \cite{adel2}. 
  For the Chern connection $\d$ attached to $q$ the curvature is given by
  \begin{equation}
  \label{orderone}
  F_{ij} =
   \frac{1}{4}\{q^{-1} ( \d_i q )q^{-1} 
  (\d_j q)-q^{-1} (\d_j q )q^{-1}(\d_i q)\}. \end{equation}
   Note that $F_{ij}$ depends only of the first (rather than the first and second) order derivatives of $q$; cf. the discussion in  \cite{dui}. 
  On the other hand 
  assume $\overline{\d}$ is a  connection commuting with  $\d$ on $A$  such that $\overline{\d}_ix=0$; then
  the $(1,1)$-curvature  of $\d$ with respect to  $\bar{\d}$ is given by 
 \begin{equation}
 \label{unknown}
  F_{i\bar{j}} = \frac{1}{2} \d_{\bar{j}}(q^{-1}\d_i q).  \end{equation}
  In particular, if $n=1$, 
 $\rho_{i\bar{j}}  =  ``\frac{1}{2} \d_i\d_{\bar{j}} \log q"$. Note that the $(1,1)$-curvature  depends on the first {\it and} second derivatives of $q$.
 
 \begin{remark} We would like to compare here our Chern connection with the Levi-Civita connection. Assume the notation above and assume, in addition,  that $n=m$. 
(More invariantly we assume that a bijection is given between the index set for the derivations $\d_i$ and the index set for the variables $x_{kl}$;
fixing such a bijection can be viewed as an analogue of Cartan's  {\it soldering} 
\cite{sternberg}, p. 351.)
One classically makes the following definition: the connection $\d$ on $G$ is {\it torsion free} if
\begin{equation}
\label{torsion}
\Gamma_{ijk}=\Gamma_{jik}.
\end{equation}
Notice the difference between the symmetry conditions \ref{symmetric} and \ref{torsion}. By the way  the system consisting of \ref{horizontal}, \ref{symmetric}, \ref{torsion}, (with $q$ given) has a solution $\Gamma_{ijk}$ if and only if $q$ satisfies the following equations:
 \begin{equation}
 \label{necessary}
 \d_i q_{jk}=\d_jq_{ik}.
 \end{equation}
 If $q$ satisfies the equations \ref{necessary} we will say $q$ is  {\it Hessian}.
 
 If $q^t=q$,  
there is a unique invariant connection $\d$ satisfying the equations
 \ref{horizontal} and \ref{torsion}; it is given by
\begin{equation}\label{windy}
\Gamma_{kij}=\frac{1}{2}\left(\d_k q_{ij}+\d_iq_{jk}-\d_jq_{ki}\right)
\end{equation}
and it is called the {\it Levi-Civita} connection attached to $q$. 
So the Levi-Civita connection and the Chern connection for a given $q=q^t$
coincide if and only if $q$ is Hessian. There is a $(1,1)$ version of this in which the Hessian condition is replaced by the {\it K\"{a}hler condition}; we will not go into this here. 

If $q^t=-q$, an invariant connection $\d$ will be called {\it symplectic} if it satisfies equations \ref{horizontal} and \ref{torsion}. A $q$ with $q^t=-q$ will be called {\it symplectic} if
\begin{equation}
\label{closed}
\d_i q_{jk}+\d_j q_{ki}+\d_k q_{ij}=0.
\end{equation}
It is trivial to see that if $\d$ is symplectic then $q$ must be symplectic.
 Also, in case $q^t=-q$,  conditions \ref{horizontal}, \ref{symmetric},  \ref{torsion} 
imply 
\begin{equation}
\label{must}
\d_k q_{ij}=0,\ \ \ \Gamma_{ijk}=0,\ \ \ \d=\d_0,\end{equation}
where we recall that $\d_0 x=0$.
The first equality  follows by combining \ref{necessary}, \ref{closed}; the second follows from  \ref{solution}; the third follows from the second. So, for $q^t=-q$,
we have that $\d$ is both symplectic and Chern if and only if $\d q_{jk}=0$ and $\d=\d_0$.

The purpose of the above review is to point out  that the Chern connection story will have an arithmetic analogue while the Levi-Civita story and the symplectic story, by themselves, will not have an arithmetic analogue in our paper. However the conjunction  of the Chern and the Levi-Civita conditions (i.e. the Hessian/K\"{a}hler condition) {\it will} have an arithmetic analogue. Also the conjunction of the Chern and symplectic conditions  will {\it also} have an arithmetic analogue.
 \end{remark}

\subsection{Arithmetic differential geometry}
Recall from \cite{char} (cf. also \cite{joyal}) that  a $p$-{\it derivation} on a $p$-torsion free ring $B$  is a set theoretic map $\d_p:B\ra B$, such that the map
$\phi_p:B\ra B$ defined by $\phi_p(b):=b^p+p\d_p (b)$ is a ring homomorphism. Note that $\phi_p$ lifts the $p$-power Frobenius on $B/pB$; we will say that $\phi_p$ is the {\it lift of Frobenius} attached to $\d_p$. Also, for any ring $B$ we denote by $B^{\widehat{p}}$ the $p$-adic completion of $B$, i.e., the projective limit of the rings $B/p^nB$.

\subsubsection{Adelic connections on $GL_n$}
Let now $A$ be the ring $\bZ[1/M,\zeta_N]$ where $M$ is some even integer and $\zeta_N$ is a primitive $N$-th root of unity, $N\geq 1$. Let $G=GL_n=Spec\ A[x,\det(x)^{-1}]$ be the general linear group scheme over $A$, where $x=(x_{kl})$ is an $n \times n$ matrix of indeterminates. So $\cO(G)=A[x\det(x)^{-1}]$. Finally let ${\mathcal V}$ be a (possibly infinite) set of primes in $\bZ$ not dividing $MN$.
We keep this notation throughout the paper. Note that for any $p$ not dividing $MN$ there is a unique $p$-derivation on $A$; indeed for any such $p$ there is a unique ring endomorphism of $A$ lifting the $p$-power Frobenius on $A/pA$; this endomorphism  sends $\zeta_N$ into $\zeta_N^p$.

\begin{definition}
An {\it adelic connection}  on $G$ is  a family $(\d_p)$, indexed by the primes $p$ in ${\mathcal V}$ where, 
for each $p$, $\d_p$ is 
a $p$-derivation on the $p$-adic completion $\cO(G)^{\widehat{p}}$ of $\cO(G)$. \end{definition}

The concept of adelic connection we just introduced, being  adelic rather than global,  should be viewed  as an analogue of the differential geometric  concept of $C^{\infty}$, rather than analytic, connection in a principal bundle. In contrast to analytic  connections which lead to trivial cohomology classes \cite{atiyah}, $C^{\infty}$ connections encode non-trivial cohomological information \cite{kobayashi} and hence our adelic connections could be viewed as encoding some  cohomology-like 
information on $Spec\ \bZ$.

Given an adelic connection $(\d_p)$ we can consider the  attached  family $(\phi_p)$ of lifts of Frobenius 
on the rings 
$\cO(G)^{\widehat{p}}$;
 we shall still denote by
$\phi_p:G^{\widehat{p}}\ra G^{\widehat{p}}$ 
the induced morphism of $p$-adic formal schemes 
where  $G^{\widehat{p}}$ is the $p$-adic completion  of $G$.

Let $T$ be the matrix $x-1$ where $1$ is the identity matrix. 

\begin{definition}
An adelic connection $\d=(\d_p)$ on $G$ {\it fixes}   $1$  if, for all $p$, $\phi_p:\cO(G)^{\widehat{p}}\ra \cO(G)^{\widehat{p}}$ sends the ideal of $1$ into itself. (Note that if this is the case then there is an  induced homomorphism $\phi_p:A^{\widehat{p}}[[T]]\ra A^{\widehat{p}}[[T]]$ for each $p$.) 
An adelic connection $\d=(\d_p)$
is said to be {\it global along $1$} if it fixes $1$ and the induced homomorphisms 
$\phi_p:A^{\widehat{p}}[[T]]\ra A^{\widehat{p}}[[T]]$
send $A[[T]]$ into $A[[T]]$. (We then still denote by $\phi_p:A[[T]]\ra A[[T]]$ the induced homomorphisms.)
\end{definition}

Morally this fits into the ideology of ``analytic continuation   between primes  introduced in  \cite{laplace}.

\subsubsection{Curvature, $3$-curvature, $(1,1)$-curvature}

For adelic connections that are global along $1$ one can introduce various concepts of curvature.

\begin{definition}
Let $\d$ be an adelic connection on $G$ that is global along $1$, with associated homomorphisms $\phi_p:A[[T]]\ra A[[T]]$. The {\it curvature} of $\d$ is the family $\varphi=(\varphi_{pp'})$, indexed by pairs of primes $p,p^{\prime}$ in $\mathcal V$, where $\varphi_{pp^{\prime}}$ are the $\bZ$-module endomorphisms
\begin{equation}
\varphi_{pp^{\prime}}:=\frac{1}{pp^{\prime}}[\phi_p, \phi_{p^{\prime}}]:A[[T]]\ra A[[T]].\end{equation}\end{definition}

Here $[\alpha ,\beta ]$ denotes the usual commutator $\alpha\circ \beta-\beta\circ \alpha$ of two $\bZ$-module endomorphisms $\alpha,\beta$ of $A[[T]]$. 
If $I\subset A[[T]]$ is the ideal generated by $T$ then for any integer $\nu\geq 1$ we denote by 
\begin{equation}
\varphi_{pp^{\prime}}^{(\nu)}:I/I^{\nu+1}\ra I/I^{\nu+1}\end{equation}
the $\bZ$-module maps induced by $\varphi_{pp'}$.

\begin{definition}
Let $\d$ be an adelic connection on $G$ that is global along $1$,
 with associated homomorphisms $\phi_p:A[[T]]\ra A[[T]]$. The {\it $3$-curvature} of $\d$ is the family $(\varphi_{pp'p''})$, indexed by triples  of primes $p,p',p''$ in $\mathcal V$, where $\varphi_{pp'p''}$ are the $\bZ$-module endomorphisms
$$\varphi_{pp'p''}:=[\phi_p,\varphi_{p'p''}]=\frac{1}{p'p''}[\phi_p,[\phi_{p'},\phi_{p''}]]:A[[T]]\ra A[[T]].$$
\end{definition}

\begin{definition}
Let $\d$ and $\overline{\d}$ be two adelic connections on $G$ that are global along $1$,
 with associated homomorphisms $\phi_p:A[[T]]\ra A[[T]]$ and $\phi_{\overline{p}}:A[[T]]\ra A[[T]]$, respectively. The {\it $(1,1)$-curvature} of $\d$ with respect to $\overline{\d}$ is the family $\varphi=(\varphi_{p\overline{p}'})$, indexed by pairs of primes $p,p^{\prime}$ in $\mathcal V$, where $\varphi_{p\overline{p}^{\prime}}$ are the $\bZ$-module endomorphisms
 \begin{equation}
\label{zanzibarr}
\varphi_{p\bar{p}^{\prime}}:=\frac{1}{p\bar{p}^{\prime}}[\phi_p, \phi_{\bar{p}^{\prime}}]:A[[T]]\ra A[[T]]\ \ \text{for $p\neq p'$ and}
\end{equation}
\begin{equation}
\label{zanzibarrr}
\varphi_{p\bar{p}}:=\frac{1}{p}[\phi_p, \phi_{\bar{p}}]:A[[T]]\ra A[[T]].\end{equation}
\end{definition}

Again, if $I\subset A[[T]]$ is the ideal generated by $T$ then for any integer $\nu\geq 1$ we denote by 
\begin{equation}
\varphi_{p\overline{p}^{\prime}}^{(\nu)}:I/I^{\nu+1}\ra I/I^{\nu+1}\end{equation}
the $\bZ$-module maps induced by $\varphi_{p\overline{p}'}$.
From now on, in this paper,  we will take $\overline{\d}=(\overline{\d}_p)=(\d_{\overline{p}})$ to be equal to the adelic connection $\d_0=(\d_{0p})$, $\d_{0p}x=0$; this adelic connection is global along $1$.

 \subsubsection{Chern connections}
Next we would like to apply the concepts above to  some basic adelic connections constructed in \cite{adel2}.
Let $q\in GL_n(A)$ with $q^t=\pm q$, where the $t$ superscript means ``transpose".
Attached to $q$ we have  maps $\cH_q:G\ra G$ and  $\cB_q:G\times G\ra G$ defined by
$$\cH_q(x)=x^tqx,\ \ \ \cB_q(x,y)=x^tqy.$$
We continue to denote by $\cH_q,\cB_q$ the maps induced on the $p$-adic completions $G^{\widehat{p}}$ and $G^{\widehat{p}}\times G^{\widehat{p}}$.  
Consider the unique adelic connection $\d_0=(\d_{0p})$ on $G$ with $\d_{0p}x=0$. Denote by $(\phi_p)$ and $(\phi_{0p})$ the families of lifts of Frobenius attached to an adelic connection $\d$ and to $\d_0$ respectively. In particular $\phi_{0p}(x)=x^{(p)}:=(x_{ij}^p)$.
The following  was  proved in \cite{adel2}:

\begin{theorem} \cite{adel2} 
 Let $q\in GL_n(A)$ be such that $q^t=\pm q$. Then there exists a unique adelic connection $\d=(\d_p)$ such that the following diagrams are commutative:
\begin{equation}
 \label{coconut}
  \begin{array}{rcl}
G^{\widehat{p}} & \stackrel{\phi_p}{\longrightarrow} & G^{\widehat{p}}\\
\cH_q  \downarrow & \  & \downarrow \cH_q \\
G^{\widehat{p}} & \stackrel{\phi_{0p}}{\longrightarrow} & G^{\widehat{p}}\\
\end{array}\ \ \ \ \ \ 
\begin{array}{rcl}
G^{\widehat{p}} & \stackrel{\phi_{0p} \times \phi_p}{\longrightarrow} & G^{\widehat{p}}\times G^{\widehat{p}}\\
\phi_p \times \phi_{0p} \downarrow & \  & \downarrow \cB_q\\
 G^{\widehat{p}}\times G^{\widehat{p}} & \stackrel{\cB_q}{\longrightarrow} & G^{\widehat{p}}\end{array}\end{equation}\end{theorem}
 
 \begin{definition}
 The adelic connection $\d$ in Theorem \ref{coconut} will be called the {\it Chern connection} (on $G=GL_n$) attached to $q$. \end{definition}

 Our aim is to show that, for certain natural $q$'s the Chern connection attached to $q$ is global along $1$ (hence has  well defined curvatures); then we will state and prove vanishing and non-vanishing results for these curvatures.

 \section{Main results}
 
 Our first result is:
 
 \begin{theorem}
 \label{globe}
Let  $q\in GL_n(A)$ be such that $q^t=\pm q$ and assume
that  all the entries of $q$ are roots of unity or $0$.  Then:

1)  The Chern connection $\d=(\d_p)$ attached to $q$ is global along $1$.

2) If $J$ is the ideal in $A[[T]]$ generated by the entries of the matrix $x^tqx-q$ and if $(\phi_p)$ are the lifts of Frobenius attached to $(\d_p)$ then $\phi_p(J)\subset J$ for all $p$.

3) Assume $n=2r$, $q=1$,  and $x=\left(\begin{array}{cc} a & b\\ c & d\end{array}\right)$ with $a,b,c,d$,  $r\times r$ matrices. Let $K$ be the ideal in $A[[T]]$ generated by the entries of the matrix $x^tx-1$ and by the entries of the matrices $a-d$ and $b+c$. Let $(\phi_p)$ be the lifts of Frobenius attached to $(\d_p)$. Then $\phi_p(K)\subset K$ for all $p$.
\end{theorem}
 
 In particular, if all the entries of $q\in G(A)$ are roots of unity or $0$, the curvature $(\varphi_{pp'})$,  $3$-curvature $(\varphi_{pp'p''})$, and  $(1,1)$-curvature $(\varphi_{p\overline{p}'})$ of $\d$ are defined. 
 
 Also we have induced $\bZ$-module endomorphisms
 $$\varphi_{pp'}^{SO}:A[[T]]/J\ra A[[T]]/J;$$
 the family $(\varphi_{pp'}^{SO})$ will be referred to as the {\it special orthogonal} curvature of the Chern connection attached to $q$.
 Of course, $A[[T]]/J$ is the ring of functions of the completion  along $1$ of the {\it special orthogonal subgroup scheme} $SO(q)$ of $G=GL_n$ (defined as the identity component of the group defined by the equations
 $x^tqx=q$.)
 
  If, in addition, $n=2r$,  $q=1$, and if $(\d_{\overline{p}})$ is the  adelic connection with $\d_{\overline{p}}x=0$ then we have induced $\bZ$-module homomorphisms
 $$\varphi_{p\overline{p}'}^{U}:A[[T]]/K\ra A[[T]]/K;$$
  the family $(\varphi_{p\overline{p}'}^{U})$ will be referred to as the {\it  unitary} $(1,1)$-curvature of the Chern connection attached to $q=1$.
 Of course, $A[[T]]/K$ is the ring of functions of the completion  along $1$ of the {\it unitary subgroup scheme} $U^c_r$ of $G=GL_{2r}$  defined by the equations
 $x^tx=1$, $a=d$, $b=-c$. (This is, of course,  an analogue of the classical complexified unitary group; see \cite{adel2}.)
 
 \begin{remark}
 The condition that the entries of a symmetric $q$ be roots of unity or $0$ can be arguably viewed as an analogue of the Chern$+$Levi-Civita$=$Hessian/K\"{a}hler condition in classical differential geometry. Also the condition that the entries of an antisymmetric $q$ be roots of unity or $0$ can be  viewed as an analogue of the Chern$+$symplectic condition in the classical case.
 \end{remark}

 Next let us say that a matrix $q\in GL_n(A)$ is {\it split} if it is one of the following:
 \begin{equation}
 \label{scorpion3}
  \left(\begin{array}{cl} 0 & 1_r\\-1_r & 0\end{array}\right),\ \ 
\left( 
\begin{array}{ll} 0 & 1_r\\1_r & 0\end{array}\right),\ \ 
\left( \begin{array}{lll} 1 & 0 & 0\\
0 & 0 & 1_r\\
0 & 1_r & 0\end{array}\right),
\end{equation}
where $n=2r,2r,2r+1$ respectively and $1_r$ is the identity $r\times r$ matrix. By Theorem \ref{globe} the Chern connection attached to any split matrix is global along $1$, hence has well defined curvature, $3$-curvature, and $(1,1)$-curvature. Our next results are about the vanishing/non-vanishing of these curvatures. 

In what follows $q\in G(A)$, $q^t=\pm q$. 

For  curvature we will prove:

 \begin{theorem}
 \label{coconutt}
 Let $q$ be  split  and let $(\varphi_{pp^{\prime}})$ be the curvature of the Chern connection on $G$ attached to $q$.
 
 1) Assume  $n\geq 4$. Then for all $p\neq p^{\prime}$ we have $\varphi^{(pp^{\prime})}_{pp^{\prime}}\neq 0$; in particular $\varphi_{pp^{\prime}}\neq 0$.
 
 2) Assume $n$ even. Then for all $p,p^{\prime}$ we have $\varphi^{(2)}_{pp^{\prime}}=0$.
 
 3) Assume $n=2$ and $q^t=-q$. Then for all $p,p^{\prime}$ we have $\varphi_{pp^{\prime}}= 0$.
 
 4) Assume $n=1$. Then for all $p,p'$ we have $\varphi_{pp'}=0$.
  \end{theorem}
 
 Morally what the above says is that for $q$ split
the curvature is non-vanishing for $n\geq 4$ (cf. assertion 1) but the non-vanishing of the curvature is somewhat subtle because it does not occur modulo $I^3$ (cf. assertion 2). For  $q=-q^t$ in case $n=2$ and for $q^t=q$ in case $n=1$ the curvature vanishes (cf. assertions 3 and 4). Note that our Theorem says nothing about whether the  curvature vanishes for $n=2,3$ in the case  $q^t=q$. Indeed our method to prove assertion 1 of the above theorem, i.e. the fact that the lifts of Frobenius $\phi_p$, corresponding to the Chern connection,  do not commute on $GL_n$ for $n\geq 4$ is to find a closed subgroup scheme $V\subset GL_n$  such that all $\phi_p$ send $V^{\widehat{p}}$ into itself and such that the restrictions of the $\phi_p$'s to $V^{\widehat{p}}$ extend to non-commuting endomorphisms of $V$.   For $n$ even  $V$ will be the ($2$-dimensional vector) group
whose points are
$$\left(
\begin{array}{cccccc}
1 & ... & 0 & 0 & ... & a\\
\  & ... & \  & \  & ... & \  \\
0 & ... & 1 & b & ... & 0\\
0 & ... & 0 & 1 & ... & 0\\
\  & ... & \  & \  & ... & \ \\
0 & ... & 0 & 0 & ... & 1
\end{array}
\right)$$
where all the diagonal entries are $1$ and all the entries except the diagonal entries and $a,b$ are $0$. For this subgroup to exist one needs  $n\geq 4$. 
The same method  will be used to prove Theorem \ref{new} below; in that case 
the role of  $V$ will be played by a subgroup of $SO(q)$ isomorphic to  the unipotent radical of $GL_3$.

For the special orthogonal  curvature we will prove:

\begin{theorem}
\label{new} 
 Let $q$ be  split  and let $(\varphi^{SO}_{pp^{\prime}})$ be the special orthogonal  curvature of the Chern connection attached to $q$.
 
1) Assume $n\geq 6$. Then for all $p\neq p'$  we have $\varphi^{SO}_{pp'}\neq 0$.

2) Assume $n=2$. Then for all $p,p'$  we have $\varphi^{SO}_{pp'}=0$.
\end{theorem}

Note that the theorem says nothing about the cases $3\leq n\leq 5$.

For the $3$-curvature we will prove:

\begin{theorem}
\label{frau101}
 Let $q=q^t$
  be  split, let $n=2r\geq 2$, and let $(
\varphi_{pp'p''})$ be the $3$-curvature of the Chern connection attached to $q$.
Then $\varphi_{pp'p''}\neq 0$ for $p'\neq p''$. 
\end{theorem}

For the $(1,1)$-curvature, we will prove:

 \begin{theorem}
 \label{mixed}
 Let $q$ be  split, let $n=2r\geq 2$   and let $(\varphi_{p\bar{p}^{\prime}})$
  be the $(1,1)$-curvature of the Chern connection  attached to $q$. 
  Then for all $p,p'$ we have $\varphi^{(1)}_{p\overline{p}'}=0$ and 
 $\varphi^{(2)}_{p\overline{p}'}\neq 0$.
   In particular $\varphi_{p\overline{p}'}\neq 0$.
 \end{theorem}
 
 Finally for the unitary  $(1,1)$-curvature we will prove:
 
    \begin{theorem}
      \label{fata}
      Let $n=2r\geq 2$ and let $\varphi^{U}_{p\overline{p}'}$ be the unitary  $(1,1)$-curvature of the Chern connection attached to attached to $q=1$. Then $\varphi_{p\overline{p}'}\neq 0$ for all $p,p'$.      \end{theorem}

\section{Proofs}

  \subsection{Explicit formula}
  Recall from \cite{adel2} that for each $p$, the $\phi_p$ in Theorem \ref{globe} 
 is induced by the unique ring homomorphism
$$\phi_p:\cO(G)^{\widehat{p}}\ra \cO(G)^{\widehat{p}},\ \ \ x\mapsto \phi_p(x)=\Phi_p,$$
prolonging $\phi_p:A^{\widehat{p}}\ra A^{\widehat{p}}$ where
 $\Phi_p$ is the 
 $n\times n$ matrix with coefficients in $\cO(G)^{\widehat{p}}$ defined by
 \begin{equation}
 \label{turandot}
\Phi_p(x)=x^{(p)}\{(x^{(p)t}\phi_p(q)x^{(p)})^{-1} (x^tqx)^{(p)}\}^{1/2}.\end{equation}
In this expression, the $t$ exponent means {\it transpose},  $a^{(p)}$ for a matrix $a$ is the matrix obtained from $a$ by raising each entry to the $p$th power, and the $1/2$ power is computed using the usual matrix series  
\begin{equation}
\label{radical}
(1+u)^{1/2}=1+\sum_{i=1}^{\infty} \left(\begin{array}{c} 1/2\\i\end{array}\right) u^i\end{equation}
 for  $u$ an $n\times n$ matrix of indeterminates. The formula is applied to
 $$u=(x^{(p)t}\phi_p(q)x^{(p)})^{-1} (x^tqx)^{(p)}-1,$$
 which has coefficients in $p\cO(G)$; this makes the series \ref{radical} $p$-adically convergent in view of the fact  that
   $\left(\begin{array}{c} 1/2\\i\end{array}\right)\in \bZ[1/2]$.

   \subsection{Proof of Theorem \ref{globe}} Since 
   the entries of $q$ are roots of unity or $0$ we have 
    $\phi_p(q)=q^{(p)}$.
 Set
 $$
 v(T)  :=  ((1+T^t)^{(p)}\phi_p(q)(1+T)^{(p)})^{-1} ((1+T^t)q(1+T))^{(p)}.$$
Now $u(T):=v(T)-1$ has entries in
$$pA^{\widehat{p}}[[T]]\cap (T)A[[T]]=p(T)A[[T]].$$
Hence the matrix
$$\Phi_p^0(T):=\Phi_p(1+T)-1=(1+T)^{(p)}\cdot v(T)^{1/2}-1,
$$
which a priori has coefficients in $A^{\widehat{p}}[[T]]$, has actually coefficients in $A[[T]]$ and has zero constant term, which ends the proof of assertion 1 of the theorem. 

To check assertion 2 note that the commutativity of the left diagram in \ref{coconut} is equivalent to the equality
$$\Phi_p(x)^t \phi_p(q)\Phi_p(x)=(x^tqx)^{(p)}.$$
So we have
$$\phi_p(x^tqx-q)=\Phi_p(x)^t \phi_p(q)\Phi_p(x)-\phi_p(q)=(x^tqx)^{(p)}-q^{(p)}.$$
But now the entries of $(x^tqx)^{(p)}-q^{(p)}$ belong to the ideal generated by 
the entries of $x^tqx-q$ in $A[x]=A[T]$.  We conclude that
$\phi_p(J)\subset J$.

To check assertion 3 note that, by the above computation 
the entries of $$\phi_p(x^tx-1)=(x^tx)^{(p)}-1$$
are in the ideal generated by  the entries of $x^tx-1$. It remains to show that 
$\Phi_p$ evaluated at the matrix $\left(\begin{array}{rr} a & b\\ - b & a\end{array}\right)$ has the form $\left(\begin{array}{rr} \alpha & \beta\\ - \beta & \alpha\end{array}\right)$, i.e., it is in the centralizer of $\left(\begin{array}{rr} 0 & 1\\ - 1 & 0\end{array}\right)$ which is clear from \ref{turandot}. 
\qed

\subsection{Some remarkable polynomials}
In some of the next arguments  certain polynomials will play a recurring  role. Let  $v,w$ be two variables and let $p, p', p''$ be odd primes. Then consider the following two polynomials (one for the upper sign and one for the lower sign):
\begin{equation}
\label{frau1}
f_p(v,w):=\frac{1}{2}(\pm(\pm v +w)^p+v^p\mp w^p)\in \bZ[1/2][v,w].
\end{equation}
Also consider the corresponding polynomials in $\bZ[1/2][u,v]$:
\begin{equation}
\label{frau10}
f_{pp'}(v,w):=f_{p'}(f_p(v,w),f_p(w,v)),
\end{equation}
\begin{equation}
\label{moloch}
g_{pp'}(v,w):=f_{pp'}(v,w)-f_{p'p}(v,w),
\end{equation}
\begin{equation}
\label{frau3}
f_{pp'p''}(v,w):=f_p(f_{p'}(f_{p''}(v,w),f_{p''}(w,v)),f_{p'}(f_{p''}(w,v),f_{p''}(v,w))),
\end{equation}
\begin{equation}
\label{frau4}
g_{pp'p''}(v,w):=f_{pp'p''}(v,w)-f_{p'pp''}(v,w)-f_{p''pp'}(v,w)+f_{p''p'p}(v,w).
\end{equation}

\begin{lemma}
\label{gargantua}
For $p\neq p'$ we have
$g_{pp'}(v,w)\neq 0$.
\end{lemma}

{\it Proof}.
We have:
\begin{equation}
\label{papa}
\begin{array}{rcl}
f_{pp'}(v,w)& := & \frac{1}{2} (v\pm w)^{p{p^{\prime}}}\\
\ &  \  & \  \\
\  & \  & +\frac{1}{2}\left(\frac{\pm (\pm v+w)^{p}+v^{p}\mp w^{p}}{2}\right)^{{p^{\prime}}}\\
\ &  \  & \  \\
\  & \  & \mp \frac{1}{2}\left(\frac{\pm (v\pm w)^{p}\mp v^{p}+w^{p}}{2}\right)^{{p^{\prime}}}.\end{array}\end{equation}
Consider the polynomial 
$f_{pp^{\prime}}(v,1)\in \bZ[1/2][v]$. We have the following congruences mod $(v^3)$ in $\bZ[1/2][v]$:
\medskip
$$\begin{array}{rcl}
f_{p{p^{\prime}}}(v,1) & \equiv & \frac{1}{2}\left(\pm 1+p{p^{\prime}}v\pm \frac{p{p^{\prime}}(p{p^{\prime}}-1)}{2}v^2\right)\\
\  & \ & \  \\
\  & \  & \mp \frac{1}{2}\left(1\pm \frac{p}{2}v+\frac{p(p-1)}{4}v^2\right)^{{p^{\prime}}}\\
\  & \  & \ \\
\  & \equiv &  \frac{1}{2}\left(\pm 1+p{p^{\prime}}v\pm \frac{p{p^{\prime}}(p{p^{\prime}}-1)}{2}v^2\right)\\
\  & \  & \  \\
\  & \  & \mp \frac{1}{2}\left(1+{p^{\prime}}(\pm \frac{p}{2} v+\frac{p(p-1)}{4}v^2)+\frac{{p^{\prime}}({p^{\prime}}-1)}{2}(\pm \frac{p}{2} v)^2\right)\\  
\  & \  & \  \\
\  & \equiv & \frac{p{p^{\prime}}}{4}v\pm \frac{p{p^{\prime}}}{16}(3p{p^{\prime}}-p-2)v^2 \ \ \ \text{mod}\ \ (v^3).
\end{array}$$
Similarly
$$f_{{p^{\prime}}p}(v,1) \equiv \frac{{p^{\prime}}p}{4}v\pm \frac{{p^{\prime}}p}{16}(3{p^{\prime}}p-{p^{\prime}}-2)v^2\ \ \ \text{mod}\ \ (v^3)$$
hence we get
\begin{equation}
\label{frau100}
g_{pp'}(v,1)\equiv \pm \frac{pp'}{16}(p'-p)v^2
\ \ \ \text{mod}\ \ \
(v^3).
\end{equation}
This ends our proof.
\qed

\medskip

In the next lemma we choose the upper sign in the definition of $f_p(v,w)$.

\begin{lemma}
\label{frau6}
$g_{pp'p''}(v,w)\neq 0$ for $p\neq p'$.
\end{lemma}

{\it Proof}.
An immediate computation gives:
$$f_{pp'p''}(v,w)=\frac{1}{2}(v+w)^{pp'p''}+\frac{1}{2}h_{pp'p''}(v,w),$$
where
$$\begin{array}{rcl}h_
{pp'p''}(v,w)& = &
\left( 
\frac{(v+w)^{p'p''}+f_{p''}(v,w)^{p'}-f_{p''}(w,v)^{p'}}{2}\right)^p\\
\  & \  & \  \\
\  & \  & -
\left( 
\frac{(v+w)^{p'p''}+f_{p''}(w,v)^{p'}-f_{p''}(v,w)^{p'}}{2}\right)^p.\end{array}$$
Next,
an easy computation as in Lemma \ref{gargantua} gives
 the following congruences mod 
$v^3$ in $\bZ[1/2][v]$:
$$h_{pp'p''}(v,w) \equiv -\left(
1+\frac{3}{4}pp'p''v+\frac{pp'p''}{32}(p'p''+2p''-12+9pp'p'')v^2\right)
\ \ \ \text{mod}\ \ \ (
v^3).
$$
Hence 
$$
\begin{array}{rcl}
2g_{pp'p''}(v,1)
 & \equiv & h_{pp'p''}(v,1)
 -h_{p'pp''}(v,1)-
 h_{p''pp'}(v,1)
 +h_{p''p'p}(v,1)\\
 \  & \  & \  \\
 \  & \equiv & -\frac{pp'p''}{32}\left(
 (p'p''+2p'')-(pp''+2p'')-(pp'+2p')+(p'p+2p)\right)v^2\\
 \  & \  & \  \\
 \  & \equiv & -\frac{pp'p''}{32}(p''-2)(p'-p)v^2
 \ \ \ \text{mod}\ \ \
 (v^3).
 \end{array}$$
which ends our proof.
\qed

\subsection{Proof of assertion 1 in Theorem \ref{coconutt}} \label{1coco}
In what follows we assume $n=2r$; the case $n=2r+1$   can be treated similarly.
Let us make the convention that the upper sign in $\pm,\mp$ refers to the  case $q^t=q$ and the lower sign refers to the case $q^t=- q$. Set $\phi_p(T)=\Phi^0_p$; so $\Phi^0_p=\Phi^0_p(T)$ is a matrix with coefficients in $A[[T]]$, hence we may evaluate $\Phi^0_p$ at various matrices. 
Let $B$ be the $r\times r$ upper corner of the $2r\times 2r$ matrix $T$.
By \ref{turandot},
  evaluating $\Phi^0_p$ at the matrix 
  $\left( \begin{array}{cc} 0 & B\\ 0 & 0\end{array}\right)$
   gives
$$
\Phi^0_p\left( \begin{array}{cc} 0 & B\\ 0 & 0\end{array}\right)
=
\left( \begin{array}{cc} 0 & F_p(B)\\ 0 & 0\end{array}\right),
$$
where
\begin{equation}
\label{musca}
F_p(B):=\frac{1}{2}(\pm(\pm B+B^t)^{(p)}+B^{(p)}\mp B^{(p)t}).\end{equation}
So $F_p(B)$ has coefficients in
 $A[B]$, hence the ring endomorphism 
 $$\phi_p:A[[B]]\ra A[[B]],\ \ B\mapsto
\phi_p(B)=
  F_p(B)$$
 sends $A[B]$ into itself. We have, of course, a commutative diagram
 \begin{equation}
 \label{did}
 \begin{array}{rcl}
 A[[B]] & \stackrel{\phi_p}{\longrightarrow} & A[[B]]\\
 u\uparrow & \  & \uparrow u\\
 A[[T]] & \stackrel{\phi_p}{\longrightarrow} & A[[T]]
 \end{array}
 \end{equation}
 where $u$ is the surjective homomorphism sending 
 $T$ into $\left( \begin{array}{cc} 0 & B\\ 0 & 0\end{array}\right)$.
 In order to conclude the proof assume  $\phi_{p}$ and $\phi_{{p^{\prime}}}$ commute on   $A[[T]]$ for some $p\neq p'$ and let us derive a contradiction.
 By \ref{did}  $\phi_{p}$ and $\phi_{{p^{\prime}}}$ commute on   $A[[B]]$. 
Clearly
we have 
 $$F_p(B_{ij})=
 F_p(B)_{ij}=\frac{1}{2}(\pm(\pm
 B_{ij}+B_{ji})^p+B_{ij}^p\mp
 B_{ji}^p)=f_p(B_{ij},B_{ji}).
 $$
 Since $\phi_{p}\phi_{{p^{\prime}}}(B)=F_{{p^{\prime}}}(F_{p}(B))$
  we have
$$\begin{array}{rcl}
\phi_{p}\phi_{{p^{\prime}}}(B)_{ij}
 & = & 
  F_{{p^{\prime}}}(F_{p}(B))_{ij}\\
 \  & \  & \  \\
 \  & = & f_{p'}(F_p(B)_{ij},F_p(B)_{ji})\\
 \  & \  & \  \\
 \  & = & f_{p'}(f_p(B_{ij},B_{ji}),f_p(B_{ji},B_{ij}))\\
 \  & \  & \  \\
 \ & = & f_{pp'}(B_{ij},B_{ji}).\end{array}$$
where $f_{pp'}$ is an in \ref{frau10}.
Choose now two indices $i\neq j$ (this is possible because $n\geq 4$). Then
we get $f_{pp'}(v,w)=
f_{p'p}(v,w)$,
i.e.,
$g_{pp'}(v,w)=0$ which contradicts
 Lemma \ref{gargantua}.
 \qed

   \subsection{Proof of Theorem \ref{frau101}}
  In the notation of subsection \ref{1coco} we have 
  $$\begin{array}{rcl}
  p'p''\varphi_{pp'p''}(B) & = &
  F_{p''}(F_{p'}(F_p(B)))-F_{p'}(F_{p''}(F_p(B)))\\
  \  & \  & \  \\
  \  & \  & 
- F_p(F_{p''}(F_{p'}(B)))+F_p(F_{p'}(F_{p''}(B))).
\end{array}$$
  Exactly as in that proof we have
  $$F_{p''}(F_{p'}(F_p(B_{ij})))=
  f_{p''p'p}(B_{ij},B_{ji})$$
  hence
  $$p'p''\varphi_{pp'p''}(B_{ij})=g_{p''p'p}(B_{ij},B_{ji})$$
  and,
  taking any two indices $i,j$ with
  $i\neq j$ (which is possible because $n\geq 4$) we conclude by
 Lemma \ref{frau6}.
  \qed
 
\subsection{Proof of Theorem \ref{new}}

 For simplicity we only treat the case $n=2r$; the case $n=2r+1$ can be reduced to the case $n=2r$. Set $x=\left(\begin{array}{cc} a & b\\ c & d\end{array}\right)=T+1$.
Let $H\subset G=GL_n$ be the subgroup scheme whose points are 
the specializations of the matrices 
$\left(\begin{array}{ll} a & 0\\ 0 & d\end{array}\right).$
By the way, 
 $H\cap SO(q)$ is the subgroup of $H$  of all matrices in $H$ such that $d=(a^t)^{-1}$. 
 The embedding $H\subset G$ corresponds to the homomorphism sending $a,d$ into themselves and $b,c$ into $0$. Set $\phi_{p}(x)=\Phi_{p}$.
 By \ref{turandot} we have
 $$\Phi_p\left(\begin{array}{ll} a & 0\\ 0 & d\end{array}\right)=
 \left(\begin{array}{ll} S_p(a,d) & 0\\ 0 & S_p(d,a)\end{array}\right),
 $$
 where
 $$S_p(a,d):=a^{(p)}\cdot \{(a^{(p)})^{-1}(d^{(p)t})^{-1}(d^t a)^{(p)}\}^{1/2}.$$
So $\phi_p:A[[T]]\ra A[[T]]$ sends the ideal of $H$ into itself. We also know that $\phi_p$ sends the ideal of $SO(q)$ into itself. Hence $\phi_p$  sends the ideal of  $H\cap SO(q)$ into itself. So
we have a commutative diagram
\begin{equation}
\label{storm}
\begin{array}{ccc}
A[[a-1]] & \stackrel{\Sigma_p}{\longrightarrow} & A[[a-1]]\\
 u \uparrow & \ & \uparrow u\\
A[[T]]/J & \stackrel{\phi_p}{\longrightarrow} & A[[T]]/J\end{array}
\end{equation}
where 
 \begin{equation}
 \label{hurricane}
u(x)=\left(\begin{array}{ll} a & 0\\ 0 & (a^t)^{-1}\end{array}\right),\ \  \Sigma_p(a):=S_p(a,(a^t)^{-1})=a^{(p)}\cdot \{(a^{-1})^{(p)}a^{(p)}\}^{-1/2}.\end{equation}
  To prove assertion 1 assume $n=2r\geq 6$ and $\varphi^{SO}_{pp'}=0$ for some $p\neq p'$ and seek a contradiction. We may assume $r=3$. In view of diagram \ref{storm} we have
 \begin{equation}
 \label{stormy}
 \Sigma_{p'}\Sigma_p(a)=\Sigma_p\Sigma_{p'}(a).
 \end{equation}
 A trivial  computation using \ref{hurricane} yields
 $$\Sigma_p\left(
 \begin{array}{lll}
 1 & u & v\\
 0 & 1 & w\\
 0 & 0 & 1\end{array}\right)=\left(
 \begin{array}{lll}
 1 & u^{p} & f_{p}(v,uw-v)\\
 0 & 1 & w^{p}\\
 0 & 0 & 1\end{array}\right);$$
 from here a trivial (but somewhat ``miraculous") computation yields:
 $$ \Sigma_p\Sigma_{p'}\left(
 \begin{array}{lll}
 1 & u & v\\
 0 & 1 & w\\
 0 & 0 & 1\end{array}\right)=\left(
 \begin{array}{lll}
 1 & u^{pp'} & f_{pp'}(v,uw-v)\\
 0 & 1 & w^{pp'}\\
 0 & 0 & 1\end{array}\right).
 $$
In the above $f_p,f_{pp'}$ are the polynomial \ref{frau1}, \ref{papa} (with the choice of the upper sign). 
 From \ref{storm} we get $f_{pp'}=f_{p'p}$ which contradicts Lemma \ref{gargantua}; this ends our proof of assertion 1.
Assertion 2 follows easily.  \qed

\subsection{Proof of assertion 2 in Theorem \ref{coconut}}\label{proofofcoco}
Let  $T=\begin{pmatrix} A & B \\ C & D \end{pmatrix}$, where $A,B,C,D$ are $r\times r$ matrices; the use of $A$ here and the use of $A$ as a name for the ring $\bZ[1/M,\zeta_N]$ should not lead to any confusion.   Denote by $(T)$ the ideal  generated by the entries of $T$. By our assumptions  $q^t=\pm q$ and $q$ is split.
The upper sign in the proof below corresponds to the case $q^t=q$ and the lower sign corresponds to the case $q^t=-q$. Set $\phi_p(T)=\Phi^0_p$ and $\varphi_{pp'}(T)=\Phi_{pp'}^0$. So
$$\Phi_{pp'}^0(T)=\frac{1}{pp'}(\Phi_{p'}^0(\Phi^0_p(T))-\Phi^0_p(\Phi^0_{p'}(T))).$$
We need to show that $\Phi_{pp'}^0\equiv 0$ mod $(T)^3$.

We proceed by computing $\Phi_p(1+T)$ modulo $(T)^3$ using formula \ref{turandot}. 
Let $1 \leq i,j,k \leq r$ and set $\underline{i} = i+r$. We have
$$
(1+T)^{(p)} = \begin{pmatrix} (1+T_{11})^p & T_{12}^p & \cdots \\ T_{21}^p & (1+T_{22})^p & \cdots \\ \vdots & \vdots &  \ddots & \end{pmatrix}
 \equiv \begin{pmatrix}E_1 & 0 \\ 0 & E_2\end{pmatrix} \,\text{mod} (T)^3,$$
where $E_1$ and $E_2$ are diagonal matrices with 
$$(E_1)_{ii} = 1+pT_{ii}+\left(\begin{array}{c} p\\ 2\end{array}\right) T_{ii}^2,\ \ \ (E_2)_{ii} = 1+pT_{\underline{i}\underline{i}}+\left(\begin{array}{c} p\\ 2\end{array}\right)T_{\underline{i}\underline{i}}^2.$$
 Therefore,
\begin{dmath*}
(1+T^t)^{(p)}q(1+T)^{(p)} \equiv \begin{pmatrix} E_1 & 0 \\ 0 & E_2 \end{pmatrix}q\begin{pmatrix} E_1 & 0 \\ 0 & E_2 \end{pmatrix} \,\text{mod} (T)^3 = \begin{pmatrix} 0 & E_1E_2 \\  \pm E_2E_1 & 0\end{pmatrix} = G,
\end{dmath*}
where the last equality defines $G$. Next, we have
$$\begin{array}{l}
((1+T^t)q(1+T))^{(p)}  =  \left(\left(\begin{array}{rcl} 1+A^t & C^t \\ B^t & 1+D^t \end{array} \right)q \begin{pmatrix}1+A & B \\ C & 1+D\end{pmatrix}\right)^{(p)}\\
 \  \\
= \left(\begin{array}{rcl} (C\pm  C^t\pm  C^tA+A^tC)^{(p)} & (1+D+A^t+A^tD\pm C^tB)^{(p)} \\ 
(\pm 1\pm  A\pm  D^t\pm  D^tA+B^tC)^{(p)} & (\pm B+B^t+B^tD\pm D^tB)^{(p)} \end{array}\right) \ \end{array}$$
Write 
\begin{equation}
\label{defofA}
A_{ij} = T_{ij}, \ \ B_{ij}=T_{i\underline{j}}, \ \ C_{ij}=T_{\underline{i}j}, \ \ D_{ij}=T_{\underline{i}\underline{j}},\end{equation}
hence
\begin{align*}
(A^tD)_{ij} &= \sum_k T_{ki}T_{\underline{k}\underline{j}}, & (D^tA)_{ij} &= \sum_k T_{\underline{k}\underline{i}}T_{kj},\\
(B^tC)_{ij} &= \sum_k T_{k\underline{i}}T_{\underline{k}j}, & (C^tB)_{ij} &= \sum_k T_{\underline{k}i}T_{k\underline{j}},
\end{align*}
and define
\begin{equation}
\label{defofQ}
Q_i=\sum_k(T_{ki}T_{\underline{k}\underline{i}}\pm T_{\underline{k}i}T_{k\underline{i}}).\end{equation}
Let $H$  be the diagonal matrix with
$$
H_{ii} =1+p(T_{ii}+T_{\underline{i}\underline{i}}+Q_i)+{p \choose 2}(T_{ii}^2+T_{\underline{i}\underline{i}}^2+2T_{ii}T_{\underline{i}\underline{i}}).
$$
Then we have
$$((1+T^t)q(1+T))^{(p)}\equiv \left(\begin{array}{rcl} 0 & H \\ \pm H & 0\end{array}\right) =F\ \ \text{mod} (T)^3,$$
where the last equality is the definition of $F$.
One trivially checks that 
$$G^{-1}F =\left(\begin{array}{cc}
E_1^{-1}E_2^{-1}H & 0\\ 0 & E_2^{-1}E_1^{-1}H\end{array}\right)
\equiv 1\ \ \ \text{mod} \ \ \ (T)^2.$$
By the binomial expansion
\begin{dmath*}(1+ (G^{-1}F-1))^{1/2} \equiv 1 + \frac{1}{2}(G^{-1}F-1) \,\ \ \ \text{mod}(T)^3.
\end{dmath*}
 So
\begin{dmath*}
\Phi_p(1+T) \equiv (1+T)^{(p)}+\frac{1}{2}(1+T)^{(p)}\frac{1}{2}(G^{-1}F-1) \\
\equiv \begin{pmatrix}E_1 & 0 \\ 0 & E_2\end{pmatrix}+\frac{1}{2}\begin{pmatrix}E_1 & 0 \\ 0 & E_2\end{pmatrix}\begin{pmatrix} E_1^{-1}E_2^{-1}H-1 & 0 \\ 0 & E_2^{-1}E_1^{-1}H-1\end{pmatrix}\\
\equiv \frac{1}{2}\begin{pmatrix} E_2^{-1}H+E_1 & 0 \\ 0 & E_1^{-1}H+E_2\end{pmatrix}\, \ \ \text{mod}\ \  (T)^3.
\end{dmath*}
Set $X=E_2^{-1}H+E_1$ and $Y=E_1^{-1}H+E_2$ and  note that $X$ and $Y$  are diagonal matrices. Then
\begin{equation}
\label{XY}
\Phi_p(1+T) \equiv \frac{1}{2}
\left( \begin{array}{cc}
X & 0 \\ 0 & Y\end{array}\right)\ \ \text{mod}\ \ (T)^3.
\end{equation}
On the other hand
\begin{equation}
\label{Xii}
\begin{array}{rcl}
X_{ii}  & \equiv  &  \big(1-pT_{\underline{i}\underline{i}}-{p \choose 2}T_{\underline{i}\underline{i}}^2+p^2T_{\underline{i}\underline{i}}^2\big)\times\\
\  & \  & \  \\
\  & \  & \times \big(1+ p(T_{ii}+T_{\underline{i}\underline{i}}+Q_i)+{p \choose 2}(T_{ii}^2+T_{\underline{i}\underline{i}}^2+2T_{ii}T_{\underline{i}\underline{i}}\big)\big)\\
\  & \  & \  \\
\  & \  & +\big(1+pT_{ii}+{p \choose 2}T_{ii}^2\big)\ \ \text{mod}\ \ (T)^3\\
\  & \  & \  \\
\  & \equiv & 2+2pT_{ii}-pT_{ii}T_{\underline{i}\underline{i}}+2{p \choose 2} T_{ii}^2+pQ_i\ \ \text{mod}\ \ (T)^3,
\end{array}\end{equation}
and similarly
\begin{equation}
\label{Yii}Y_{ii}\equiv 2+2pT_{\underline{i}\underline{i}}-pT_{ii}T_{\underline{i}\underline{i}}+2{p \choose 2} T_{\underline{i}\underline{i}}^2+pQ_i\ \ \text{mod}\ \ (T)^3.\end{equation}
By the way equations \ref{XY}, \ref{Xii}, and \ref{Yii} provide a formula for $\Phi_p(1+T)$ mod $(T)^3$ and, in particular, one gets  that $\Phi_p(1+T)$ is congruent modulo $(T)^3$ to a diagonal matrix.

Now, in order to compute $\Phi_{pp'}^0$ mod $(T)^3$ set
$Z_p=\Phi_p(1+T)=\Phi_p^0(T)+1$.
Since $Z_p$ is congruent mod $(T)^3$ to a diagonal matrix and since for any $n\times n$ diagonal matrix $\Delta$ and any $n\times n$ matrix $\Gamma$ we have $(\Gamma \Delta)^{(p^{\prime})}=\Gamma^{(p^{\prime})}\Delta^{(p^{\prime})}$, the following holds:
\begin{dmath*}
\Phi_{p^{\prime}}^0(\Phi_p^0(T)) = Z_p^{(p^{\prime})}\left[\left(Z_p^{t(p^{\prime})}qZ_p^{(p^{\prime})}\right)^{-1}(Z_p^tqZ_p)^{(p^{\prime})}\right]^{\frac{1}{2}}-1\\
 \equiv Z_p^{(p^{\prime})}\left[\left((Z_p^{t}qZ_p)^{(p^{\prime})}\right)^{-1}(Z_p^tqZ_p)^{(p^{\prime})}\right]^{\frac{1}{2}}-1\ \ \text{mod}\ \ (T)^3\\
 = Z_p^{(p^{\prime})}-1.
\end{dmath*}
By \ref{Xii} we get
\begin{equation}
\label{see}(\frac{1}{2}X_{ii})^{p^{\prime}}\equiv 1+pp^{\prime}T_{ii}-\frac{pp^{\prime}}{2}T_{ii}T_{\underline{i}\underline{i}}+\frac{pp^{\prime}}{2}Q_i+\frac{pp^{\prime}(pp^{\prime}-1)}{2}T_{ii}^2\ \ \text{mod}\ \ (T)^3,\end{equation}
and similarly, by \ref{Yii}, we get
\begin{equation}
\label{seee}
(\frac{1}{2}Y_{ii})^{p^{\prime}}\equiv 1+pp^{\prime}T_{\underline{i}\underline{i}}-\frac{pp^{\prime}}{2}T_{ii}T_{\underline{i}\underline{i}}+\frac{pp^{\prime}}{2}Q_i+\frac{pp^{\prime}(pp^{\prime}-1)}{2}T_{\underline{i}\underline{i}}^2\ \ \text{mod}\ \ (T)^3.\end{equation}
Consequently
$$Z_p^{(p^{\prime})}\equiv Z_{p^{\prime}}^{(p)}\ \ \text{mod}\ \ (T)^3,$$
hence $\Phi^0_{pp^{\prime}}\equiv 0$ mod $(T)^3$.
\qed

\begin{remark}
\label{noninvertible} 
As already noticed in the proof above, $\Phi^0_p$ is congruent mod $(T)^3$ to a diagonal matrix.
 Note, on the other hand,  that $\Phi^0_p$ is {\it not} congruent mod $(T)^4$ to a diagonal matrix (e.g. for $p=3$). By the way, modulo $(T)^2$, the situation is  simpler; we have
$$\Phi^0_p\equiv \text{diag}(pT_{11},...,pT_{nn})\ \ \ \text{mod}\ \ (T)^2.$$ 
 \end{remark}

\subsection{Proof of Theorem \ref{mixed}}
Recall that by our assumptions $n=2r\geq 2$ and  $q^t=\pm q$ is split. We let $\Phi_{p\overline{p}'ij}^0:=\varphi_{p\overline{p}'}(T_{ij})$ be the entries of the matrix
$\Phi^0_{p\overline{p}'}:=\varphi_{p\overline{p}'}(T)$ where $i,j=1,...,n=2r$. As before, we set $\underline{i}=i+r$ for $1\leq i\leq r$. Also we use the notation in \ref{defofA} and \ref{defofQ}. 
Our theorem will be proved if we prove the following congruences:

\medskip

1) $\Phi^0_{p\overline{p}'ii}\equiv \Phi^0_{p\overline{p}'\underline{i}\underline{i}}\equiv \frac{1}{2}(Q_i-T_{ii}T_{\underline{i}\underline{i}})$ mod $(T)^3$ for $1\leq i\leq r$ and $p\neq p'$.

\medskip

2) $\Phi^0_{p\overline{p}ii}\equiv \Phi^0_{p\overline{p}\underline{i}\underline{i}}\equiv \frac{p}{2}(Q_i-T_{ii}T_{\underline{i}\underline{i}})$ mod $(T)^3$ for $1\leq i\leq r$.

\medskip

3) $\Phi^0_{p\overline{p}'ij}\equiv 0$ mod $(T)^3$ for $1\leq i,j\leq n=2r$,  $i\neq j$.

\medskip

\noindent To check these we use the notation in subsection \ref{proofofcoco}. Set  
$$Z_{\overline{p}'}:=\Phi_{\overline{p}'}(1+T)=\Phi_{\overline{p}'}^0(T)+1=(1+T)^{(p')}.$$
 Then exactly as in 
subsection \ref{proofofcoco} we have 
$$\Phi^0_{\overline{p}'}(\Phi^0_p(T))=Z_p^{(p')}-1,$$
$$\Phi^0_p(\Phi^0_{\overline{p}'}(T))=Z_{\overline{p}'}^{(p)}-1.$$
Now $Z_p^{(p')}$ mod $(T)^3$ was computed in \ref{see} and \ref{seee} while  $Z_{\overline{p}'}^{(p)}$ mod $(T)^3$ is, of course, trivial to compute; then assertions 1, 2, 3 above  follow. 
\qed

\medskip

Since assertion 4 in Theorem \ref{coconut} trivially follows from the formula \ref{turandot} all we need to conclude the proof of Theorem \ref{coconut}  is:

\subsection{Proof of assertion 3 in Theorem \ref{coconutt}}
 We need to show that $\phi_p$ and $\phi_{p^{\prime}}$ commute on $A[[T]]$ for $q^t=-q$, $q$ split.  In this case formula \ref{turandot} yields
$$\Phi_p(x)=\lambda_p(x)\cdot x^{(p)}, \ \ \ \ \lambda_p(x)=\left(
\frac{\det(x^{(p)})}{(\det(x))^p}
\right)^{-1/2}.$$
Since $x^{(pp^{\prime})}=x^{(p)(p^{\prime})}=x^{(p^{\prime})(p)}$ and
 $$\phi_p(T)=\lambda_p(1+T)\cdot (1+T)^{(p)}-1$$
 we get:
$$
\begin{array}{rcl}
\phi_{p}(\phi_{p'}(T)) & = & \lambda_{p^{\prime}}(\lambda_p(x)\cdot x^{(p)})\cdot (\lambda_p(x)\cdot x^{(p)})^{(p^{\prime})}-1\\
\  & \  & \  \\
\   & = & \lambda_{p^{\prime}}(x^{(p)})\cdot (\lambda_p(x))^{p^{\prime}}\cdot x^{(pp^{\prime})}-1\\
\  & \  & \  \\
\  & = & \left(
\frac{\det(x^{(pp^{\prime})})}{\det(x^{(p)})^{p^{\prime}}}
\right)^{-1/2}\cdot \left(
\frac{\det(x^{(p)})^{p^{\prime}}}{(\det(x))^{pp^{\prime}}}
\right)^{-1/2}\cdot x^{(pp^{\prime})}-1\\
\  & \  & \  \\
\  & = & \left(
\frac{\det(x^{(pp^{\prime})})}{(\det(x))^{pp^{\prime}}}
\right)^{-1/2}\cdot x^{(pp^{\prime})}-1\\
\  & \ & \  \\
\  & = & \phi_{p'}(\phi_{p}(T)), 
\end{array}
$$
which proved the desired commutation.
\qed

\medskip

To conclude the proofs of our theorems we need:

\medskip

 {\it Proof of Theorem \ref{fata}}. We may assume $M=2$ and $N=1$ so our ground ring is $\bZ[1/M,\zeta_N]=\bZ[1/2]$. Again, write
$T=\begin{pmatrix} A & B \\ C & D \end{pmatrix}$ and $x=\left(\begin{array}{cc} a & b\\ c & d\end{array}\right)$ so $a=1+A$, $b=B$, $c=C$, $d=1+D$. Denote by  $GL^c_r$ the subgroup  of $GL_{2r}$ defined by the equations $a=d$, $b=-c$. So 
$U^c_r$ is defined in $GL^c_r$ by the equations $aa^t+bb^t=1$, $ba^t=ab^t$, equivalently by the equations:
$$F:=A+A^t+AA^t+BB^t,\ \ G:=(B+BA^t)-(B^t+AB^t).$$
Similarly  $GL_1^c=Spec\ \bZ[1/2][\alpha,\beta,\frac{1}{(1+\alpha)^2+\beta^2}]$, with $\alpha,\beta$ two variables, and 
$U^c_1$ is defined in $GL_1^c$
 by the equation
$$f:=2\alpha+\alpha^2+\beta^2.$$
Consider the embedding of $GL_1^c$ into $GL_r^c$ given by the  homomorphism between their coordinate rings
$$u:\left(\begin{array}{rr} A & B\\ -B & A\end{array}\right)\mapsto \left(\begin{array}{rr} \alpha \cdot 1_r& \beta\cdot 1_r\\ -\beta\cdot 1_r & \alpha\cdot 1_r\end{array}\right).$$
The ideal of the image of $GL_1^c$ via this embedding is easily seen to be sent into itself by the lifts of Frobenius $\phi_p$ on $(GL_r^c)^{\widehat{p}}$; hence the ideal of $U^c_1$ is sent into itself by the lifts of Frobenius $\phi_p$ on $(U^c_r)^{\widehat{p}}$. So we have 
a commutative diagram
 \begin{equation}
 \label{diddd}
 \begin{array}{ccc}
 \bZ[1/2][[\alpha,\beta]]/(f) & \stackrel{\phi_p}{\longrightarrow} & \bZ[1/2][[\alpha,\beta]]/(f)\\
 u\uparrow & \  & \uparrow u\\
 \bZ[1/2][[A,B]]/(F,G) & \stackrel{\phi_p}{\longrightarrow} & \bZ[1/2][[A,B]]/(F,G)\\
  v\uparrow & \  & \uparrow v\\
 \bZ[1/2][[T]] & \stackrel{\phi_p}{\longrightarrow} & \bZ[1/2][[T]]
 \end{array}
 \end{equation}
 where  $v$ sends $A$ into $A$, $B$ into $B$, $C$ into $- B$ and $D$ into $A$.  An easy computation yields
 $$\phi_p\left(\begin{array}{rr} \alpha & \beta\\ -\beta & \alpha\end{array}\right)=\left(
 \begin{array}{rr} \alpha_p & \beta_p\\ - \beta_p & \alpha_p\end{array}\right)$$
 where
 $$\alpha_p=(1+\alpha)^pK_p(\alpha,\beta)-1,\ \ \ \beta_p=\beta^pK_p(\alpha,\beta),$$
 $$K_p(\alpha,\beta):=\left(\frac{((1+\alpha)^2+\beta^2)^p}{(1+\alpha)^{2p}+\beta^{2p}}\right)^{1/2}.$$
On the other hand we have
$$\phi_{\overline{p}'}
\left(\begin{array}{rr} \alpha & \beta\\ -\beta & \alpha\end{array}\right)=
\left(\begin{array}{rr} (1+\alpha)^{p'}-1 & \beta^{p'}\\ -\beta^{p'} & (1+\alpha)^{p'}\end{array}\right).$$
Assume $\varphi_{p\overline{p}'}^{U}=0$ and seek a contradiction.
 The equality $\varphi_{p\overline{p}'}^{U}=0$ implies that $\phi_p$ and $\phi_{\overline{p}'}$ commute on 
 $\bZ[1/2][[\alpha,\beta]]/(f)$
 which immediately implies
 $$K_p(\alpha,\beta)^{p'}=K_p((1+\alpha)^{p'}-1,\beta^{p'}).$$
 Squaring the latter equation and  using the relation
 $\beta^2=-2\alpha-\alpha^2$, 
   we get the following equality of polynomials in $\alpha$:
    $$(1+\alpha)^{2pp'}-(2\alpha+\alpha^2)^{pp'}=((1+\alpha)^{2p}-(2\alpha+\alpha^2)^{p})^{p'}
   \cdot ((1+\alpha)^{2p'}-(2\alpha+\alpha^2)^{p'})^{p}.
   $$
   Picking out the coefficients of $\alpha$ we get 
   $2pp'= 4pp'$, a contradiction; this ends the proof.
 \qed
 
 \bigskip

\section{The case of one prime}

So far there was no restriction on the cardinality of our ``index set" of primes ${\mathcal V}$. Interestingly, restricting to the case when ${\mathcal V}$ has one element, is still interesting and actually allows one to introduce $(1,1)$-curvature
for connections that are not necessarily global along the identity. We give in what follows some details.

As in the rest of the paper we let $A=\bZ[1/M,\zeta_N]$ with $M$ even.
Consider the set ${\mathcal V}=\{p\}$ consisting of one prime $p$ only, such that $p$ does not divide $MN$. 
In this case an adelic connection is simply a $p$-derivation $\d_p:\cO(G)^{\widehat{p}}\ra \cO(G)^{\widehat{p}}$.
Now let us consider two adelic connections $\d_p$ and $\overline{\d}_p=\d_{\overline{p}}$ (not necessarily global along $1$). To these one can attach the {\it $(1,1)$-curvature}:
  \begin{equation}
  \label{zanzibar}
  \varphi_{p\overline{p}}:=\frac{1}{p}[\phi_p,\phi_{\overline{p}}]:\cO(G)^{\widehat{p}}\ra \cO(G)^{\widehat{p}}.\end{equation} 
 In case $\d$ and $\overline{\d}$ are global along $1$, 
  \ref{zanzibar}  is compatible (induces)  \ref{zanzibarrr}.    Note that  $\Phi_{p\overline{p}}:=\phi_p(x)$ is an $n\times n$ matrix with coefficients in $\cO(G)=A[x,\det(x)^{-1}]^{\widehat{p}}$ hence the value, $\Phi_{p\overline{p}}(1)$, of
$\Phi_{p\overline{p}}$ at $x=1$ is an $n\times n$ matrix with entries in $A^{\widehat{p}}$. Under the assumptions above we have:

  \begin{theorem}
   \label{sunny}
   Let  $q\in GL_n(A)$, $q^t=\pm q$.
   Let $\d=(\d_p)$ be  the Chern connection  attached to $q$ and let $\overline{\d}=(\d_{\overline{p}})$ be the adelic connection with $\d_{\overline{p}}x=0$; here both adelic connections are indexed by the set $\{p\}$ with one prime $p$.  Let $\varphi_{p\overline{p}}$  be the $(1,1)$-curvature  of  $\d$  with respect to  $\overline{\d}$ and let $\Phi_{p\overline{p}}:=\varphi_{p\overline{p}}(x)$.  Then: 
   
   1) $\Phi_{p\overline{p}}(1)=-\d_{\overline{p}}\left(\left(1+p(q^{(p)})^{-1}\d_pq
\right)^{-1/2}\right).$

2) If $n=1$ then $\Phi_{p\overline{p}}=\Phi_{p\overline{p}}(1)\cdot x^{p^2}$.

3) $\Phi_{p\overline{p}}(1)= 0$ if and only if each entry of $q$ is  either $0$ or a root of unity. Moreover,  if $n=1$, then $\varphi_{p\overline{p}}=0$ if and only if $q$ is a root of unity.
   \end{theorem}
   
   {\it Proof}.
   The formulae in assertions 1 and 2   follows from a direct computation using \ref{turandot}. Assertion 3 follows from assertions 1 and 2 plus the following two facts: i)   if $\d_pa=0$ or $\d_{\overline{p}}a=0$ for some $a\in A$ then $a$ is a root of unity or $0$; and ii) if $1+pa$ is a root of unity or $0$ for some $a\in A$ then $a=0$.
   \qed
   
\begin{remark}
Note the analogy between the formula in assertion 1 of Theorem \ref{sunny} and formula \ref{unknown}.
\end{remark}

\end{document}